\documentclass[12pt]{article}
\usepackage[utf8]{inputenc}
\usepackage{amssymb}
\usepackage{amsmath}
\usepackage[left=2cm,right=2cm, top=2cm,bottom=2cm]{geometry}

\newcommand{\cH}{\mathcal{H}}

\newcommand{\cM}{{\mathcal{M}}}

\newcommand{\fU}{{\mathfrak{U}}}
\newcommand{\tr}{\operatorname{tr}}

\newcommand{\RR}{{\mathbb{R}}}

\newcommand{\CC}{{\mathbb{C}}}

\newcommand{\bC}{\boldsymbol{C}}

\newcommand{\ld}{\ldots}

\title{Some parabolic equations for measures\\and Gaussian semigroups}
\author{O.\,E.~Galkin, S.\,Yu.~Galkina}
\date{November 2020}
\begin{document}

\maketitle

\centerline{\it Department of Fundamental Mathematics,}
\centerline{\it Laboratory of Dynamical Systems and Applications}

\centerline{\it National Research University Higher School of Economics} 
\bigskip

\textbf{Abstract.} This short communication (preprint) is devoted to mathematical study of evolution equations that are important for mathematical physics and quantum theory; we present new explicit formulas for solutions of these equations and discuss their properties. The results are given without proofs but the proofs will appear in the longer text which is now under preparation.

In this paper, infinite-dimensional generalizations of the Euclidean analogue of the 
Schr\"odinger equation for anharmonic oscillator are considered in the class of measures. 
The Cauchy problem for these equations is solved. In particular cases, 
explicit formulas for fundamental solutions are obtained, which are a generalization 
of the Mehler formula, and the uniqueness of the solution with certain properties is proved. 
An analogue of the Ornstein-Uhlenbeck measure is constructed. 
The definition of Gaussian semigroups is given and their connection with the considered parabolic equations is described.

\textbf{Keywords:} linear evolution equation, Schr\"odinger equation, parabolic equation, equations for measures, Cauchy problem, exact solution, Gaussian semigroup 

\textbf{MSC2010 codes:} \textbf{
28C20, 
35C15, 
35K15, 
35Q40, 
37B55
}

\textbf{Email:} ogalkin@hse.ru, sgalkina@hse.ru

\section{Introduction}

Quantum mechanics \cite{Moretti} is one of the basic tools of modern science and technology. Quantum computers, small objects in molecular biology, processes involving nano-particles that are important e.g.\ for constructing new materials with useful properties -- all obey the rules of quantum mechanics. Laws of quantum mechanics in its modern understanding are written on the language of functional analysis. In the paper we heavily employ this language and can recommend the following books as a source of terminology and facts that are used without references: \cite{BS-RFA, HF}.

One of the main questions answered by mathematical theories is predicting the future of the system (particle, several particles, airplane, robot, chemical reactor, field, quantum computer etc) based on the information about the laws of evolution (i.e. rules of changing in time) and current state of the system considered. For example, Newton's laws of classical mechanics help to predict position and velocity of a classical particle based on the knowledge of the acting forces and also initial position and velocity of the particle. In this simple example the second Newton's law is the evolution equation (i.e. equation setting the rules of evolution), all possible values of position and velocity is the space of states of the system, and initial position and velocity is the state of the system in the initial moment. The evolution equation together with the initial state are called the initial-value problem, or Cauchy problem, and the solution of this problem is a function that gives the state of the system for all points of non-negative timeline. General theory of evolution equations is closely related with the branch of functional analysis called operator semigroups theory (see for example~\cite{Remizov2018}); in the paper we use notions of this theory without references, but all the facts and terminology we use can be found in standard textbooks \cite{HF, EN}.

In quantum mechanics the main evolution equation  is the Schr\"odinger equation (see\cite{LanLif1977, GlJa1987}). This is an equation of the form $ih\cdot\psi'_t=H\psi$, where $\psi$ is the state of a quantum system and $H$ is the quantum Hamiltonian -- linear self-adjoint operator in Hilbert space. The Schr\"odinger equation is in some sense the quantum analogue the Newton's equation in classical mechanics, the Hilbert space describes states of the quantum system, and the Hamiltonian describes the process of quantum evolution. The equation $\psi'_t=-H\psi$ is called the Euclidean analogue of the Schr\"odinger equation; in the case considered in the paper it is a second order parabolic partial differential equation that is connected with both classic and quantum evolution. In the paper the role of $\psi$ is played by the Borel measure \cite{Kuo1975, Bogachev1998, BS-RFA}, theory of such measures is a part of infinite dimensional analysis \cite{IDA, Dapr, BS-TVS}.

We call \emph{a $n$-dimensional anharmonic oscillator} a system consisting of $n$ one-dimensional interacting oscillators.
The Schr\"odinger equation for the wave function 
$\psi(t;x_1,\ldots,x_n)$ of a quantum $n$-dimensional anharmonic oscillator has the form (\cite{LanLif1977}):
\begin{equation} \label{eqUrShrGO1}
ih\cdot\psi'_t = -\frac12\sum_{j=1}^{n} \frac{\psi''_{x_jx_j}}{m_j} + \frac12\sum_{j=1}^{n} k_jx_j^2\cdot\psi + V\cdot\psi,
\end{equation}
where $m_1,\ldots,m_n$ are the masses of one-dimensional oscillators, $k_1,\ldots,k_n$ are the stiffness of their springs.
This equation can be rewritten in a more general way.
To do this, we will define two diagonal matrices: the matrix $B_n$ with diagonal elements $1/(hm_j)$, $j=1,\ld,n$, and the matrix $C_n$ with diagonal elements $k_j/h$, $j=1,\ld,n$.
In addition, we denote by $\psi''_{xx}$ the Hesse matrix with elements 
$\psi''_{x_jx_k}$, $j,k=1,\ld,n$. Then the equation~\eqref{eqUrShrGO1}
can be written as:
\begin{equation} \label{eqUrShrGO2}
\psi'_t = i\Big[\frac12\tr(B_n\psi''_{xx}) - \frac12(C_nx,x)\cdot\psi - V\cdot\psi\Big].
\end{equation}
Removing the imaginary unit on the right side, we get the Euclidean analog of the equations \eqref{eqUrShrGO1} and \eqref{eqUrShrGO2}: 
\begin{equation} \label{eqUrShrGO3}
\psi'_t = \frac12\tr(B_n\psi''_{xx}) - \frac12(C_nx,x)\cdot\psi - V\cdot\psi.
\end{equation}
Next, we turn to the infinite-dimensional analog of this equation.
In the infinite-dimensional case, there is, in particular, the following fundamental difference from the $n$-dimensional case:
in an infinite-dimensional Hilbert space, there is no analog of the Lebesgue measure (see~\cite{Weil1951}).
The presence of the Lebesgue measure is important because in quantum mechanics, the square of the modulus of the wave function is equal to the density of the coordinate distribution relative to this measure.
The absence of its analog in the infinite-dimensional case leads to the idea of using not only the wave function $\psi$, 
but also the "wave" measure $\mu$ to describe a quantum system. 
In this case, the probability of the falling of system coordinates into the set $A$, is equal to the integral $\int_A\psi\,d\mu$.

Thus, we come to the following problem: solve the analog of the equation~\eqref{eqUrShrGO3} with respect to measures defined on an infinite-dimensional Hilbert space.
More precisely, we solve equations of the next form in our work:
\begin{equation} \label{eqFullEq}
\mu^{\,\prime}_t=
 \frac 12\operatorname{tr}(B\mu^{\,\prime\prime})
-\Bigl[(D\mu^{\,\prime},\cdot)+\operatorname{tr}D\cdot\mu\Bigr]
- \frac 12 (C\cdot,\cdot)\mu + V\cdot\mu.
\end{equation}
To determine the differential operators included in the right-hand side of the equation, as well as to solve the equation itself, we use the Fourier transform method (see~\cite{AvSm1972}). This method seems to us the most convenient.
We give more precise definitions in the next section.

\section{Preliminaries and problem setting}

Let ${\mathcal{H}}$ be a separable Hilbert space,
$B_{\mathcal{H}}$ is $\sigma$--algebra of its Borel subsets,
${\mathcal{M}}({\mathcal{H}})$ is the set of Borel measures on $B_{\mathcal{H}}$,
${\mathfrak{U}}_{\mathcal{H}}$ is the algebra of cylindrical subsets in ${\mathcal{H}}$,
${\mathcal{M}}_c({\mathcal{H}})$ is the family of cylindrical measures on ${\mathfrak{U}}_{\mathcal{H}}$.
Next let $B$, $C$ and $D$ be bounded linear operators on ${\mathcal{H}}$, at that $B$, $C$ are symmetric and nonnegative, and $B$ is nuclear.

\medskip
{\bf Definition 1.}
We will call the function
$\widetilde\mu(y) = \int_{\mathcal{H}} e^{i(x, y)} d\mu(x)$, where $y\in {\mathcal{H}}$,
Fourier transform of the measure $\mu\in{\mathcal{M}}_c ({\mathcal{H}})$.

\medskip
{\bf Definition 2.}
If $\mu\in {\mathcal{M}}({\mathcal{H}})$, then let us denote 
through $\operatorname{tr}(B\mu^{\,\prime\prime})$
such a measure $\lambda \in {\mathcal{M}}_c({\mathcal{H}})$, that 
$\widetilde\lambda (\varphi)=-(B\varphi,\varphi) \widetilde\mu (\varphi)$\
(of course, if it exists), and through
$[(D\mu^{\,\prime},\cdot)+\operatorname{tr}D\cdot\mu]$
--- such $\nu \in {\mathcal{M}}_c ({\mathcal{H}})$, that 
$\widetilde\nu(\varphi)= - \widetilde\mu^{\,\prime}_{D\varphi}(\varphi)$
(here the derivative is understood in the sense of Gato).

\medskip
{\bf Definition 3.}
The family of measures $\{\mu(t)\}_{t>0}$ is called {\it differentiable} 
({\it weakly differentiable\/}) if for all $A\in\fU_\cH$ and $t>0$ there exists 
$\dfrac{d}{dt} \bigl([\mu(t)](A)\bigr)$ (respectively, if for all functions
$f_c$ lying in the class $\bC_{bc}(\cH,\RR)$ of continuous bounded 
cylindrical functions from $\cH$ to $\RR$, 
there exists $\dfrac{d}{dt} \int_{\cH} f_c (x)\, d[\mu(t)](dx)$\ ).

\medskip
{\bf Remark 1.}
It is known that if the family $\mu (t)$ is differentiable, then the set function\newline
$\dfrac{d}{dt} \mu(t)\colon A\mapsto \dfrac{d}{dt}\bigl([\mu(t)](A)\bigr)$ is a cylindrical measure
(see~\cite{DanfSv1958}).

\medskip
{\bf Definition 4.}
Let ${\mathcal{L}}: {\mathcal{M}}({\mathcal{H}})\to {\mathcal{M}}_c ({\mathcal{H}})$ 
be some linear operator
with the domain $D_{\mathcal{L}}$.
The family of measures $\{\mu(t)\}_{t>0}$, lying in $D_{\mathcal{L}}$,
is called {\it weak solution of the equation
$\mu^{\,\prime}_t={\mathcal{L}}\mu$
with initial conditions
$\mu_0\in {\mathcal{M}}({\mathcal{H}})$},
if the next two conditions are met:

1) $\frac{\textstyle d}{\textstyle dt\rule{0pt}{1.8ex}} 
\int_{\mathcal{H}} f_c(x)\,d[\mu(t)](x)=\int_{\mathcal{H}} f_c(x) \, d[{\mathcal{L}}\mu(t)](x)$
for any continuous cylindrical bounded function $f_c$ on ${\mathcal{H}}$;

2) $\lim_{t\to + 0} \int_{\mathcal{H}} f(x)\,d[\mu(t)](x)=\int_{\mathcal{H}} f \, d\mu_0$
for any continuous bounded function $f$ on ${\mathcal{H}}$.

\medskip
{\bf Definition 5.}
We call {\it the fundamental solution} of such an equation,
the family of its solutions is $G_x(t)$, $t\geqslant 0$, depending on the parameter $x\in {\mathcal{H}}$,
with initial conditions $G_x(0)=\delta_x$.

\medskip
First let consider the following evolution equation:
\begin{equation} \label{eq1.2}
\mu^{\,\prime}_t =
 \frac 12\operatorname{tr}(B\mu^{\,\prime\prime})
-\Bigl[(D\mu^{\,\prime},\cdot)+\operatorname{tr}D\cdot\mu\Bigr]
- \frac 12 (C\cdot,\cdot)\mu + \alpha\mu.
\end{equation}
We will look for its fundamental solution in the form of
\begin{equation} \label{eq1.3}
G_x(t)=s(t)\exp\Bigl\{-\frac 12(P(t)x,x)\Bigr\}
\Gamma(Q(t),R(t)x),
\end{equation}
where
\newline
1) function $s: \RR_+\to \CC$ is continuous;
$P$, $Q$ and $R$ are continuous mappings from $\RR_+$
to $L(\cH,\cH)$, with $P\geqslant 0$, and $P$, $Q$ are symmetric;
\newline
2) $\Gamma(Q(t), R(t) x)$ is Gaussian measure with  the correlation
operator $Q(t)$ and the mean $R(t)x$.
\newline
3) $s (0)=1$, $P (0)=0$, $Q (0)=0$, $R (0)=I$
(here $I$ is a unit operator).

\medskip
{\bf Definition 6.}
Family of measures
$G_x(t)$, $t\geqslant 0$, of the form~(2) we call {\it a Gaussian semigroup},
if it satisfies the above conditions 1) -- 3),
and also has the semigroup property: \newline
$ \int_{\cH} G (t) (y) G(s) (x) (dy)=G(t+s) (x)$,
for any $t,s\geqslant 0$.

\section{Results}
\subsection{Explicit formulas for solutions of evolution equation}

{\bf Theorem 1.}
{\sl If the weak fundamental solution of~\eqref {eq1.2} has the form~\eqref{eq1.3} and has properties 1), 2) and 3), then the functions $P$, $Q$, $R$ and $s$
satisfy a system of equations and initial conditions
\begin{equation}
\left\{
\begin{aligned}
s\,' &=-\frac12s\cdot \tr(CQ)+\alpha s\\
P\,' &=R^{\,*}CR\\
Q\,' &=B-QCQ+D^{\,*}Q+QD\\
R\,' &=-QCR+D^{\,*}R
\end{aligned}
\right.\qquad\left\{
\begin{aligned}
s(0) &=1\\
P(0) &=0\\
Q(0) &=0\\
R(0) &=I.
\end{aligned} \right.
\end{equation}
}

\medskip
{\bf Theorem 2.}
{\sl If $P$, $Q$, $R$ and $s$ are the solution of~(3), then the family of measures of the form~\eqref{eq1.3} is the (strong) fundamental solution of the equation~\eqref{eq1.2}
and is a Gaussian semigroup. 
Also, if $\mu_0\in \cM(H)$ and $ \int_{\cH} |x|^2 \, d\mu < \infty$, then the family
$\mu(t)= \int_{\cH} G_x(t)\mu_0(dx)$, $t>0$, is the solution of~\eqref{eq1.2} with the initial condition $\mu_0$.}

\medskip
The validity of the following two theorems can be easily verified by direct substitution:

\medskip
{\bf Theorem 3.}
{\sl If
$C=0$, then
$$
\begin{cases}
s(t)=e^{\alpha t},&\cr
P(t)\equiv 0,&\cr
Q(t)=\int_0^t \exp\{D^{\,*}s\}B \exp\{Ds\}\,ds, &\cr
R(t)=\exp\{D^{\,*}t\} &
\end{cases}
$$
}

\medskip
{\bf Theorem 4.}
{\sl If
$D=0$, then
$$
\begin{cases}
s(t)=e^{\alpha t}\cdot \det^{-1/2}(\cosh(t\sqrt{CB})), &\cr
P(t)= C(\tanh( t\sqrt{CB}))/\sqrt{CB}, &\cr
Q(t)=\bigl(\tanh(t\sqrt{CB})/\sqrt{CB}\bigr) B \text{ and} &\cr
R(t)=\cosh^{-1}( t\sqrt{CB}) &
\end{cases}
$$
}

\medskip
{\bf Remark 2.}
Even functions of $\sqrt{CB}$ are defined using Maclaurin series, which, due to parity, will contain only integer powers of the operator $CB$.

{\bf Theorem 5.}
{\sl Let $G_x(t)$ be a Gaussian semigroup. Then there are operators $B$, $C$, and $D$ belonging to $L(\cH,\cH)$, with $B=B^{\,*}$, $C=C^{\,*}$, and nuclear $B \geqslant 0$, and a number $\alpha$ such that for the equation they define of the form~(2), $G_x(t)$ will be the fundamental solution.}

\subsection{One statement about uniqueness}

Using some ideas from~\cite{Smol1979}, we prove that the solution of (1) is unique. The result is stated in the following theorem.

\medskip
{\bf Theorem 6.}
{\sl Let $C=0$,
$\eta\in \cM(H)$. 
Then, if $ \int_{\cH} |x|^2 \, d\eta <  \infty$, then there is at most one solution $\{\mu(t)\}_{t>0}$ of the equation~\eqref{eq1.2} with the initial condition $\eta$ satisfying the condition
\newline
$ \sup_{0 < t\leqslant T} \int_{\cH} |y|^2 \|\mu(t)\|(dy) <  \infty$ for any $T > 0$ 
(here $\|\mu(t)\|$ is variation of measure $\mu(t)$).
In~addition, a solution satisfying this condition exists.
}

\subsection{Construction of conditional generalized Ornstein-Uhlenbeck measure}

Let's denote by $G_x(t)$ the fundamental solution of the equation
$$
\mu^{\,\prime}_t =
 \frac 12\operatorname{tr}(B\mu^{\,\prime\prime})
-\Bigl[(D\mu^{\,\prime},\cdot)+\operatorname{tr}D\cdot\mu\Bigr] - \frac 12 (C\cdot,\cdot)\mu.
$$
Let the symbol $\boldsymbol{C}^T_{x, A}$ means, 
for fixed $T>0$, $x\in {\mathcal{H}}$ and $A\in B_{\mathcal{H}}$,
the space of all continuous functions $f$ on the segment $[0, T]$, 
taking values in ${\mathcal{H}}$, such that $f(0)=x$ and $f(T)\in A$.
For any set of points $0<t_1<\ldots<t_n<T$ and sets $A_1,\ldots,A_n\in B_{\mathcal{H}}$,
the subsets
$I^{t_1,\cdots, t_n}_{A_1,\cdots, A_n}
 =\Bigl\{\ f \bigm| f(t_i)\in A_i,\ i=1,\ldots, n\Bigr\}$
form a semiring $K$ in $\boldsymbol{C}^T_{x, A}$.
Let's set the measure $U^T_{x, A}$ on it by the equality
$$
\begin{aligned}
&U^T_{x,A} (I^{t_1,\cdots,t_n}_{A_1,\cdots,A_n})
 =\int_{A_1}[G_x(t_1)](dy_1)
\int_{A_2}[G_{y_1}(t_2-t_1)](dy_2)\cdot\ldots\cdot\\
&\cdot
\int_{A_n}[G_{y_{n-1}}(t_n-t_{n-1})](dy_n)
[G_{y_n}(T-t_n)](A).
\end{aligned}
$$

\medskip
{\bf Theorem 7.}
{\sl 
The measure  $U^T_{x,A}$ on $K\subset \boldsymbol{C}^T_{x,A}$
has unique countably-additive continuation to Borel measure on $\boldsymbol{C}^T_{x,A}$.
Moreover, his continuation $U^T_{x,\cH}$ is a gaussian measure.
}

\medskip
{\bf Definition 7.} We will call 
{\it a conditional generalized Ornstein --- Uhlenbeck measure}
the measure $U^T_{x,A}$, described in theorem~7 (see also~\cite{Gaveau1981}).

\subsection{Existence of a solution of the equation with potential}

Here we use just constructed measure $U^T_{x,A}$ on the set of functions $\boldsymbol{C}^T_{x,A}$ 
(such functions can be called \emph{trajectories} in the Hilbert space ${\mathcal{H}}$) 
to represent solutions of the evolutionary equation of type~\eqref{eqFullEq} as an integral (see for example~\cite{Mazzucchi2009}).

{\bf Theorem 8.}\ {\sl
Let $V:{\mathbb R}_+\times {\mathcal{H}} \to {\mathbb{C}}$ be a function, 
continuous by the totality of variables, 
and there are functions $C_1, C_2 \in L_{1, loc}({\mathbb R}_+)$
and number $0\leqslant r\leqslant 2$,
such that for any $t\in{\mathbb R}$, $x\in {\mathcal{H}}$, the next inequalities are satisfied:
$$
|V(t,x)|\leqslant C_1(t)\exp\bigl\{\overline o(\|x\|)\bigr\}\quad \text{and}\quad
\operatorname{Re} V(t,x)\leqslant C_2(t) \overline o(\|x\|^r).
$$
Then the equation
$\nu\,'_t= \frac12\operatorname{tr}(B \nu'') - \bigl[(D\nu',\cdot)+\operatorname{tr} D\,\nu\bigr] + V\nu$
has a weak fundamental solution $G^V_x(t)$ of the form
$$
[G^V_x(t)](A)= \int_{\boldsymbol{C}^T_{x,A}}
\exp \Bigl\{\int_0^tV(s,q(s))\,ds\Bigr\} U^t_{x,A}(dq),
$$
where $A\in B_{\mathcal{H}}$.
Also, if $\nu_0\in{\mathcal{M}} ({\mathcal{H}})$ and
$\int_{\mathcal{H}} \exp\{\|x\|^r\}\nu_0(dx) < \infty$,
then the formula $[\nu(t)](A)=\int_{\mathcal{H}} [G^V_x(t)](A)\nu_0(dx)$, $A\in B_{\mathcal{H}}$,
sets the weak solution of this equation with the initial condition $\nu_0$.}

\section*{Acknowledgements}

Authors are thankful to I.D.~Remizov for attention to the research and for fruitful discussion of the issues raised in it. 
The publication was prepared within the framework of the Academic Fund Program at HSE University in 2020–2021 (grant No.20-04-022, project Evolution semigroups and their applications)  and within the framework of  the Russian Academic Excellence Project ``5-100''.




\end{document}